\documentclass[]{article}
\def\be{\mathbf{e}}   
\def\bbf{\mathbf{f}}   
\def\bff{\mathbf{f}}   
\def\bi{\mathbf{i}}   
\def\bj{\mathbf{j}}   
\def\bk{\mathbf{k}}   
\def\bv{\mathbf{v}}   
\def\bF{\mathbf{F}}   
\def\bH{\mathbf{H}}   
\def\bJ{\mathbf{J}}   

\def\C{\mathbb{C}} 
\def\G{\mathbb{G}} 
\def\H{\mathbb{H}}  

\def\R{\mathbb{R}}  

\def\no{\noindent}

\def\beq{\begin{equation}}
\def\eeq{\end{equation}}

\usepackage{array}
\usepackage{tabularx}
\usepackage{amsfonts}
\usepackage{amssymb}
\usepackage{tikz-cd}
\usetikzlibrary{matrix,arrows,decorations.pathmorphing}
\usepackage{graphicx}   
\def\bpm{\begin{pmatrix}}
\def\epm{\end{pmatrix}}
\makeindex            
\begin{document}

\title{Geometrization of the Real Number System}
\author{Garret Sobczyk
\\ Universidad de las Am\'ericas-Puebla
 \\ Departamento de F\'isico-Matem\'aticas
\\72820 Puebla, Pue., M\'exico
\\ http://www.garretstar.com}
\maketitle
\begin{abstract} Geometric number systems, obtained by
extending the real number system to include new {\it anticommuting} square roots of
$\pm 1$, provide a royal road to higher mathematics by largely sidestepping the tedious languages of tensor analysis and category theory. 
The well known consistency of real and complex matrix algebras, together with Cartan-Bott periodicity, firmly establishes the consistency of these geometric number systems, often referred to as Clifford algebras. The {\it geometrization} of the real number system is the culmination of the thousands of years of human effort at developing ever more sophisticated and encompassing number systems underlying scientific progress and advanced technology in the 21st Century. Complex geometric algebras are also considered.
 
\smallskip
\no {\em AMS Subject Classification:} 15A63, 15A66, 81R05, 81R25
\smallskip

\no {\em Keywords:} Cartan-Bott periodicity, Clifford algebra, complex numbers, geometric algebra, Hurwitz-Radon numbers, quaternions, real numbers, spinors. 
\end{abstract}
\newtheorem{thm1}{Theorem}
\newtheorem{thm2}[thm1]{Theorem}
\section{Introduction}

The concept of number has played a decisive role in the ebb and flow of civilizations
across centuries. Each more advanced civilization has made its singular contributions
to the further development, starting with the natural ``counting numbers" of ancient
peoples, to the quest of the Pythagoreans that (rational) numbers are
everything, to the heroic development of the ``imaginary" numbers for solving cubic and quartic polynomials, and their further generalization to Hamilton's quaternions \cite{TD67}. I maintain that
the culmination of this development is the {\it geometrization} 
of the number concept:

\begin{quote}{\bf Axiom:}   The real number system can be geometrically extended to include new,
anti-commutative square roots of $\pm 1$, each new such root representing
the direction of a unit vector along orthogonal coordinate axes of a Euclidean or
pseudo-Euclidean space $ \R^{p,q}$, where $p$ and $q$ are the number of new square roots of $+1$ and $-1$, respectively. 
\end{quote}

 The resulting real geometric algebra, denoted by
\[  \G_{p,q} := \R(\be_1, \ldots, \be_p, \bff_1, \ldots, \bff_q), \]
has dimension $2^{p+q}$ over the real numbers $\R$, and is said to be {\it universal} since no further relations between the new square roots are assumed.
Since $\G_{p,q}$  is an associative ring, it is natural
to consider matrix rings over  $\G_{p,q}$. Indeed, the elements of a geometric algebra provide a natural geometric basis for matrix algebra,  and taken together form an integrated framework which is more
powerful than either when considered alone \cite{SNF}, \cite{S08}.
   Also considered are complex geometric algebras and their corresponding complex matrix algebras.  
 The antecedents of our geometric algebras can be found in the works of W. K. Clifford \cite{wkc}, H. Grassmann \cite{grass1862}, and W. Hamilton \cite{wrh}.

The importance of geometric algebras in physics was first recognized by 
Brauer and Weyl \cite{b/w1935}, and Cartan \cite{cartan1966}, particularly in connection with the concept of $2$- and $4$-component spinors at the heart of the newly minted {\it quantum mechanics} \cite{sources1967}. The concept of a spinor arises naturally in Clifford algebras, and  much work has been carried out by mathematicians and physicists since that time \cite{dss1992}. The fundamental importance of Clifford algebras has also been recognized in the computer science and engineering communities, as well as in efforts to develop the mathematics of quantum computers \cite{HD02}.

\section{The geometric algebra $\G_{p,q}$} 

The associative geometric algebra 
\beq  \G_{p,q} = \R(\be_1, \ldots, \be_p, \bff_1, \ldots, \bff_q) .  \label{geoalg} \eeq
 Each $\be_i$ is a new square root of $+1$ for $1\le i \le p$, and each $\bff_j$ is
a new square root of $-1$ for $1 \le j \le q$.
Thus,
  \[  \be_i^2 = 1 = -\bff_i^2,   \]
and the new square roots of $\pm 1$ are  pairwise anti-commutative.

 By the {\it standard basis} of $\G_{p,q}$, we mean
\[  \G_{p,q} = span_\R\{ 1,{\cal V}^1, \ldots, {\cal V}^{n} \}, \]
where each ${\cal V}^k$ consists of $\pmatrix{n \cr k}$ products of $k$ distinct basis unit vectors  selected from the $n=p+q$ orthogonal unit vectors. Each element of ${\cal V}^k\equiv\G_{p,q}^k $ is a {\it $k$-vector}, and taken together span the homogeneous subspace ${\cal V}^k$ of $k$-vectors in $\G_{p,q}$.
As a graded real linear space, the universal geometric algebra $\G_{p,q}$ has 
\[  2^{n} = \pmatrix{n \cr 0} +\pmatrix{n \cr 1} + \cdots + \pmatrix{n \cr n}  \]
dimensions as previously mentioned. 

As an important example, the case $p=1=q$ gives the geometric algebra $\G_{1,1}$. We have
\beq \G_{1,1} = span_\R\{ 1, \be, \bff, \be \bff \} =span_\R \pmatrix{1 \cr \be}\pmatrix{1 & \bff}= span_\R \pmatrix{1 & \bff \cr \be & \be \bff}, \label{stdbasis11} \eeq
where $\be^2=1=-\bff^2$ and $\be\bff=-\bff \be$.
The equation (\ref{stdbasis11}) specifies the geometric algebra $\G_{1,1}$ in terms of
its {\it standard basis}. The last two terms on the right express the standard basis as a matrix of basis elements defined by the matrix product of the column $\pmatrix{1 \cr \be}$ with the row $\pmatrix{1 & \bff}$. In terms of the standard basis, any
element $g \in \G_{1,1}$ is given by
\beq  g = a_0 + a_1 \be + a_2 \bff + a_3 \be \bff =\pmatrix{1 & \be}\pmatrix{a_0 & a_2 \cr a_1 & a_3} \pmatrix{1 \cr \bff},\label{gstdbasis11} \eeq
where $a_\mu \in \R$ for $0 \le \mu \le 3$. In (\ref{gstdbasis11}), we see how the matrix notation
can be utilized to express $g$ in terms of the standard basis of $\G_{1,1}$.

Since the argument used to obtain the corresponding isomorphic matrix algebra $M_2(\R)$ of $\G_{1,1}$ will be generalized to higher dimensional geometric algebras, it is worthwhile to explore the relationship between the matrix algebra $M_2(\R)$ and the geometric algebra $\G_{1,1}$ in detail before proceeding further.
The unit bivector $u:=\be \bff$ satisfies 
\[  u^2 =(\be \bff)(\be \bff)=-(\be \bff)(\bff \be)=-\be^2 \bff^2 = 1,\]
and is used to define the {\it mutually annililating idempotents} $u_\pm := \frac{1}{2}(1\pm u)$.
Whereas (\ref{gstdbasis11}) utilizes the matrix product to express a general element $g\in \G_{1,1}$ in terms of the standard basis (\ref{stdbasis11}), it is not particularly useful in relating matrix multiplication to the geometric product. For this we use the {\it spectral basis} of
$\G_{1,1}$, defined by
\beq \G_{1,1}=span_\R\{1,\be, \bff, \be \bff \} =\pmatrix{1 \cr \be}u_+ \pmatrix{1 & \be} = \pmatrix{u_+ & \be u_-
\cr \be u_+ & u_-}. \label{specbasis11} \eeq

Noting that $u u_\pm =\pm u_\pm$, respectively, and $\be u_+ = u_- \be$, we get the relation
\[  \pmatrix{1 & \be}u_+ \pmatrix{1 \cr \be} =u_+ + \be u_+ \be = u_+ + u_- = 1. \]
Using this relation, we find for $g$ given in (\ref{gstdbasis11}),
\[  g=  \pmatrix{1 & \be}u_+ \pmatrix{1 \cr \be}g \pmatrix{1 & \be}u_+ \pmatrix{1 \cr \be} \]
\[ = \pmatrix{1 & \be}u_+ \pmatrix{g  & g \be \cr  \be g & \be g \be}u_+  \pmatrix{1 \cr \be}  \]
\[=  \pmatrix{1 & \be}u_+ \pmatrix{a_0+a_3 & a_1-a_2 \cr a_1 +a_2 & a_0-a_3} \pmatrix{1 \cr \be}= \pmatrix{1 & \be}u_+ [g] \pmatrix{1 \cr \be} , \]
where $[g]:= \pmatrix{a_0+a_3 & a_1-a_2 \cr a_1 +a_2 & a_0-a_3}\in M_2(\R)$ is the matrix of $g$ with respect to the spectral basis (\ref{specbasis11}).

The unique property of the spectral basis (\ref{specbasis11}) is that the geometric multiplication of two elements 
$g,h \in \G_{1,1}$
 corresponds to the corresponding matrix product of their respective
matrices $[g]$ and $[h]$, \cite{S0}. That is
\[  gh =  \pmatrix{1 & \be}u_+ [g] \pmatrix{1 \cr \be} \pmatrix{1 & \be}u_+ [h ] \pmatrix{1 \cr \be} \]
\[=  \pmatrix{1 & \be} [g]u_+ \pmatrix{1 & \be \cr \be & 1}u_+ [h]  \pmatrix{1 \cr \be} = \pmatrix{1 & \be}u_+ [g][h] \pmatrix{1 \cr \be}, \]
which establishes the algebra 
isomorphism $\G_{1,1} \widetilde =M_2(\R)$.

\section{Geometric algebra building blocks}

   We have shown in the previous section that   $\G_{1,1}\widetilde = M_2(\R)$. It is natural to make the even more basic identifications ${\bf R} \equiv \G_{0,0} \equiv M_1(\R)$ and 
   ${\bf C} \equiv \G_{0,1} \equiv M_1(\C)$.
   It is also common to identify Hamilton's quaternions $\H $  with the geometric algebra
   $\G_{0,2}$, that is $\bH \equiv{\H} \equiv \G_{0,2}$, and the
   {\it hyperbolic numbers} $\,^2{\bf R} \equiv \G_{1,0}$.
   The hyperbolic numbers, a subalgebra of $\G_{1,1}$, are naturally identified with the {\it diagonal} matrices of $M_2(\R)$, that is
   $\,^2{\bf R} \widetilde = M_2^D(\R)$. It is interesting to note that whereas the complex numbers $\C $ were first used in the $14^{\rm th}$ Century hunt for solutions to the cubic and quartic polynomial equations, the hyperbolic numbers $^2{\bf R} $ could have served the same purpose \cite{S1}.
   
   Let us examine the quaternions $\H$ more
   closely. We have
   \beq \G_{0,2}=\R(\bi,\bj)= span_\R \{1,\bi,\bj,\bk\}, \label{qcase} \eeq
  where $\bk:=\bi \bj$. It is easily checked that the elements
  $\bi, \bj, \bk $ satisfy the usual rules for quaternion multiplication. The quaternions also arise naturally as the {\it even subalgebra} $\G_{3,0}^+ $ of
  the geometric algebra 
  \[ \G_{3,0}:= span_\R \{1,\be_1, \be_2,\be_3,
  \be_{12} ,\be_{13} ,\be_{23} ,\be_{123} \} ,\] 
  where $\be_{123}:=\be_1 \be_2 \be_3$, 
  \beq  \bi :=\be_{23}  ,\ \bj :=\be_{31}, \quad {\rm and} \quad \bk := \bi \bj= \be_{21}.  \label{defquaternion} \eeq  
  
  The geometric algebra $\G_{3}\equiv \G _{3,0}$ has the spectral basis 
  \beq \G_3=\pmatrix{1 \cr \be_1}u_+ \pmatrix{1 & \be_1}= \pmatrix{u_+ & \be_1u_- \cr \be_1u_+ & u_-}_{\C}  \label{specbasis3} \eeq 
  over the formally complex numbers 
  $\C :=span_\R\{1,i \}$ where here $u:=\be_3$, $i:=\be_{123}$, and the mutually annihilating idempotents $u_\pm := \frac{1}{2}(1\pm \be_3)$. Any element $g \in \G_3$ can be written
  \[  g= \pmatrix{1 & \be_1}u_+[g]\pmatrix{1 \cr \be_1},\]
  where $[g]$ is the {\it matrix} of $g$ with respect to the spectral basis (\ref{specbasis3}). Using this spectral basis,
  \[ [\be_1]=\pmatrix{0 & 1 \cr 1 & 0},\  [\be_2]=
  \pmatrix{0 &-i \cr i & 0 },\ [\be_3] = \pmatrix{1 & 0 \cr 0 & -1},\]
  known as the {\it Pauli matrices} \cite{paulimat}. 
  From the Pauli matrices for the standard basis vectors of $\G_3$, using (\ref{defquaternion}), we calculate the corresponding
  matrix representations of the basis quaternions
  $ [\bi ]=[\be_2][\be_3]$, $[\bj]  =[\be_3][\be_1]$, and $[\bk] = [\be_2][\be_1]$,  
  
  \[  [\bi ]=\pmatrix{0 & i \cr i & 0}, \  [\bj ]=\pmatrix{0 & 1 \cr -1 & 0}, \  [\bk ]=\pmatrix{-i & 0 \cr 0 & i}.  \]
  We denote this correspondence by ${\bf H}=M_2^q(\C) $.
  
  There is one final building block that we need,
  the {\it double quaternions} 
  \[ ^2\bH:=\G_{0,3} =\R(\bff_1, \bff_2, \bff_3 ). \]
  The standard basis of $\G_{0,3}$ is
  \[\G_{0,3}=span_\R \{ 1, \bff_1, \bff_2, \bff_3,\bff_{12} , \bff_{13} , \bff_{23} , \bJ \}, \]
  for $\bJ:=\bff_{123}$, and where the quaternions are again identified as elements of the
  even subalgebra  
  \[ \G_{0,3}^+:=span_\R \{1,\bi := \bff_{23},\ \bj :=\bff_{13},\ \bk := \bff_{21}\}. \]
  A general element $Q \in \G_{0,3}$ has the form
  \[ Q = q_1+q_2\bJ ,\] 
  where $q_1:=a_0 + a_1 \bi +a_2 \bj +a_3 \bk$ and 
  $q_2:=b_0 + b_1 \bi +b_2 \bj +b_3 \bk$ for
  $a_\mu,b_\mu \in \R $.
  Noting that ${\bJ}^2=1$, we define the annihilating idempotents 
  $\bJ_\pm := \frac{1}{2}(1\pm \bJ)$. Then,
  since ${\bf J}{\bf J}_+={\bf J}_+$ and
  ${\bf J}{\bf J}_-=-{\bf J}_-$, 
  \[  Q = Q( \bJ_+ + \bJ_-)=(q_1+q_2)\bJ_+ + (q_1-q_2)\bJ_-= q_+ \bJ_++ q_- \bJ_-,  \]
  for $q_+ , q_- \in \G_{0,3}^+$ where $q_+:=q_1+q_2$ and $q_-:=q_1-q_2$. It follows that the matrix $[Q]$ of $Q$
  can be written as a {\it block diagonal} $D$-matrix,
  \[ [Q]= \pmatrix{[q_+] & [0] \cr [0] & [q_-]}\in M_2^{D}({\bf H})=\, ^2\bH, \]
  for the $D$-blocks $[q_+],[q_-] \in M_2^q(\C)$.
  
  Now that we have identified the corresponding isomorphic real and complex matrix algebras
  \[ M_1(\R), M_1(\C),M_2^D(\R), M_2^q(\C),
  M_2^{D}(\bf H) \]
   of the base geometric algebra building blocks
  \[ \G_{0,0},\G_{0,1},\G_{1,0},\G_{0,2},\G_{0,3}, \] respectively,   we proceed to develop the recursive relationship that will allow us to express any geometric algebra
  $\G_{p,q}$ as an algebra of matrices over the real or complex numbers of the building blocks. By doing so we have fully justified our geometrization of the real number system as per our Axiom given at the beginning.

\section{Classification of Geometric algebras}
The purpose of this section is to 
establish general recursive relationships between geometric algebras and matrices. We begin with 

\begin{thm1} The geometric algebra $\G_{p+1,q+1}$ is
algebraically isomorphic to the matrix geometric algebra $M_2(\G_{p,q})$, i.e., 
\beq \G_{p+1,q+1} \widetilde= M_2(\G_{p,q})  \label{bott1}  \eeq
where $p\ge 0$ and $q \ge 0$.
\end{thm1}

\no {\bf Proof:}

 Recall that 
\[\G_{p,q}=\R(\be_1, \ldots, \be_p, \bff_1, \ldots, \bff_q) \ \ {\rm and} \ \  
 \G_{p+1,q+1}=\R(\be_1, \ldots, \be_p,\be,  \bff_1, \ldots, \bff_q, \bff ). \]
 Analogous to (\ref{gstdbasis11}), any element $G\in \G_{p+1,q+1}$ can be expressed in the form
\beq  G = g_0 +  g_1\be + g_2 \bff  + g_3 \be \bff,  \label{Gspecbasiskk} \eeq
where $g_\mu \in \G_{p,q}$ for $0 \le \mu \le 3$. Applying the spectral basis (\ref{specbasis11}) to $\G_{p+1,q+1}$, and again noting that 
\[  \pmatrix{1 & \be}u_+ \pmatrix{1 \cr \be} =u_+ + \be u_+ \be = u_+ + u_- = 1, \]
we calculate

\[  G =   \pmatrix{1 & \be}u_+ \pmatrix{1 \cr \be }G \pmatrix{1 & \be}u_+ \pmatrix{1 \cr \be }\] 

\[  = \pmatrix{1 & \be}u_+ \pmatrix{G  & G \be \cr  \be  G & \be G \be}u_+  \pmatrix{1 \cr \be}   \]
\[ = \pmatrix{1 & \be}u_+ \pmatrix{g_0 + g_3  & g_1-g_2 \cr g_1^- +g_2^-  & g_0^- -g_3^- } \pmatrix{1 \cr \be}=  \pmatrix{1 & \be}u_+[G] \pmatrix{1 \cr \be},  \]
where 
\[ [G]:=  \pmatrix{g_0 + g_3  & g_1-g_2 \cr g_1^- +g_2^-  & g_0^- -g_3^- } \in M_2(\G_{p,q}), \]
and $g^-:=\be g \be$ is the operation of {\it geometric inversion} in $\G_{p,q}$ obtained by replacing all vectors in $g$ by their negatives.   

$\hfill \ensuremath{\square} $

There is another very useful relationship between
geometric algebras. Noting that for $p\ge 0$ and $q \ge 0$,
\beq \R(\be, \be_1, \ldots,\be_p,\bff_1, \ldots , \bff_q) =\R(\be,\be \bff_1, \ldots ,\be  \bff_q, \be \be_1, \ldots,\be \be_p),   \label{evensubalgebra} \eeq
it follows that 
\beq \G_{p+1,q} = \G_{q+1,p}.   \label{bott2} \eeq
From (\ref{evensubalgebra}) and (\ref{bott2}), we also have that each element 
$G \in \G_{p+1,q}$ can be written in the form
\[ G = g_1 + \be g_2   \ \ {\rm for} \ \ g_1,g_2 \in \G_{q,p}=\R(\be \bff_1, \ldots ,\be  \bff_q, \be \be_1, \ldots,\be \be_p),  \]
which shows that 
 \beq  \G_{p+1,q}^+ = \R(\be \bff_1, \ldots ,\be  \bff_q, \be \be_1, \ldots,\be \be_p) =  \G_{q,p}.  \label{evensubalgebra2} \eeq
 Similarly, since
 \[ \G_{p,q+1}= \R( \bff\be_1  , \ldots,\bff\be_p, \bff , \bff\bbf_1  , \ldots , \bff \bff_q )=\R(\be_1, \ldots ,\be_p ,\bff, \bff_1, \ldots,\bff_q),     \] 
\[ G = g_1 + \bbf g_2   \ \ {\rm for} \ \ g_1,g_2 \in \G_{p,q}=\R(\bff \be_1  , \ldots,\bff\be_p , \bff \bbf_1 , \ldots ,\bff \bff_q ),  \]
it follows that
  \beq  \G_{p,q+1}^+ = \R(\be_1 \bbf, \ldots ,\be_p \bbf , \bbf_1 \bbf , \ldots,\bbf_q \bbf )  =  \G_{p,q},  \label{evensubalgebra3} \eeq
 relating the even sub-algebras of $ \G_{p+1,q}$ and $ \G_{q+1,p}$ to
 $ \G_{q,p}$ and $ \G_{p,q}$,   respectively.

The equalities of geometric algebras in (\ref{bott2}), (\ref{evensubalgebra2}) and 
(\ref{evensubalgebra3})
are used in a loose sense, in-so-far as the elements that
are identified as generating basis vectors in $\G_{q+1,p}$ are a mixture of a vector and bivectors in
$\G_{p+1,q}$. More properly, we say that the algebras are {\it algebraically isomorphic} and write
   \[ \G_{p+1,q} \widetilde = \G_{q+1,p}. \]

In addition to the basic relationships (\ref{bott1}) and (\ref{bott2}), it is easy to establish that
\beq \G_{p+4,q}=\G_{p,q+4}. \label{bott3} \eeq 
To see this, write
\[ \G_{p+4,q}= \R(\be_a,\be_b,\be_c,\be_d,\be_1, \ldots , \be_p, \bff_1, \ldots, \bff_q  ),    \]
and
\[ \G_{p,q+4}= \R(\be_1, \ldots , \be_p,\bff^\prime_a,\bff^\prime_b,\bff^\prime_c,\bff^\prime_d, \bff_1, \ldots, \bff_q  ),    \]
where 
\[ \bff^\prime_a:= \be_b \be_c \be_d,\ \bff^\prime_b :=  \be_c \be_d \be_a ,\ \bff^\prime_c:= \be_d \be_a \be_b,\ \bff^\prime_d := \be_a \be_b \be_c .  \]
For $s\in \{a,b,c,d\}$, the $\bff^\prime_{s}\in \G_{p+4,q}$ are {\it anticommuting trivectors} which also anticommute with the vector generators of 
 $\G_{p,q}$. They also serve as anticommuting vector generators of $ \G_{p,q+4}$, the product of any distinct three of them producing $\pm$ a basis vector in $\G_{p+4,q}$.
 For example, 
 \[ \bff^\prime_a\bff^\prime_b\bff^\prime_c =-\be_ d \in \G_{p+4,q}^1.  \]
 
 As a final recursive relationship between geometric algebras and ring matrices, we have
 
 \begin{thm2} If $p-q=1(mod\ 4)$ for $p \ge 0$, $q\ge 0$, then $\G_{p+k,q}=\G_{p,q+k}$ for any integer $k$ such that $p+k \ge 0$ and $q+k \ge 0$. 
 \end{thm2}
 
 \no {\bf Proof:} Since $p-q=1(mod\ 4)$, it follows that for some $s$, $p=q+1+4s\ge 0$. Using (\ref{bott1}), (\ref{bott2}), and (\ref{bott3}), we have for $s \ge 0$,
 \[ \G_{p+k,q} = \G_{q+1+k+4s,q}=Mat(2^q,\G_{k+1+4s,0}).  \] 
But also for $s\ge 0$,
 \[ \G_{p,q+k}=\G_{q+1+4s,q+k}=Mat(2^q,\G_{1+4s,k}),\]
and
\[ \G_{1+4s,k}=\G_{1,k+4s}=\G_{k+1+4s,0}. \]
If $s <0$, then
$q=p-1-4s$ and the argument can be repeated substituting in for $q$ rather than for $p$.
 
$\hfill \ensuremath{\square} $

From the Cartan periodicity relations
 
   \beq \G_{p+8,q}= \G_{p+4,q+4}=M_{2^4}(\G_{p,q}) =
   \G_{p,q+8}, \label{cartan-period8} \eeq
 and for $p-q=1(mod\,4)$,
 \[ \G_{p+k,q}=\G_{p,q+k} \]
  established in Theorem 2, the famous Classification Table for geometric
algebras for $n=p+q$ follows; the rows are numbered by $p+q$ where $0 \le p+q \le 7$, and the columns by $p-q$ where $-7 \le p-q\le 7$. We have included an extra row in the Table 1 to give the sign of the square of the pseudoscalar $i^2=(-1)^{\frac{(p-q)(p-q-1)}{2}}$ for each of the geometric algebras $\G_{p,q}$; this is important in relationship to the {\it Hurwitz-Radon numbers} discussed in the next section. In order to fit Table 1 to the page, we have left out the parenthesis around the matrix arguments. This also helps to bring out the relationship of Table 1 to the Budinich/Trautman ``Clifford Clock" \cite{budtraut}, given in Table 2.

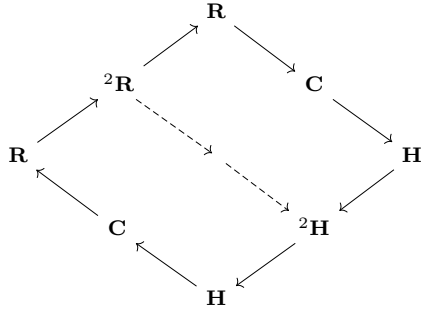
\begin{table}
	\centering
	\caption{Classification of real geometric algebras $\G_{p,q}$.}
	\label{tab:table1}

{\footnotesize 
  \[ \begin{tabular*}{\columnwidth}[t]{@{}c@{}c@{}c@{}c@{}c@{}c@{}c@{}c@{}c@{}c@{}c@{}c@{}@{}c@{}c@{}c@{}}
 7 & 6  & 5&4 &3&2&1&0 &-1&-2& -3&-4&  -5&-6 &-7 \\
0$\quad \quad\quad \quad$    && & &  & & & {\bf R} & & & & & &&  \\
 1$\quad \quad\quad \quad$ &&&&&& $^2{\bf R}$ &&  {\bf C} &&&&&& \\
  2$\quad \quad\quad \quad$  &&&&& $M_2$\bf R && $M_2${\bf R} && 
   {\bf H} & & &&& \\
3$\quad \quad\quad \quad$  &&&& $M_2${\bf C} && $M_2$$^2{\bf R}$ &&  $M_2${\bf C} && $^2{\bf H}$ &&&& \\
4$\quad \quad\quad \quad$ &&& $M_2${\bf H}&& $M_4${\bf R}&&$M_4${\bf R}&&$M_2${\bf H}&& $M_2${\bf H}\\
5$\quad \quad\quad \quad$ &&$ M_2$$^2{\bf H}$ && $M_4${\bf C}&&$M_4$$^2{\bf R}$ && $M_4${\bf C}&& $M_2$$^2{\bf H}$&& $M_4${\bf C} \\
6$\quad \quad\quad \quad$ & $M_4${\bf H} && $M_4${\bf H} && $M_8${\bf R}&&$M_8${\bf R} && $M_4${\bf H}&& $M_4${\bf H}&& $M_8${\bf R} \\
 7$\quad \ $  $M_8$\bf C && $M_4$$^2{\bf H}$ && $M_8${\bf C} && $M_8$$^2{\bf R}$&&$M_8${\bf C} && $M_4$$^2{\bf H}$&& $M_8${\bf C}&& $M_8$$^2{\bf R}$ \\
 $i^2 \quad \quad \  -$ & $-$ &$+$ &$ +$ & $-$ & $-$ & $+$ & $+$ &$ -$ & $-$ &$ +$ &$ +$ &$ -$ &$ -$ &$ +$ \\
 \end{tabular*}  \] }
\end{table}

\begin{table}[h!]
	\centering
	\caption{Budinich/Trautman Clifford Clock.}
	\label{tab:table2} 
	\bigskip
{\footnotesize
\begin{tikzcd}
	&& \bf R  \arrow{rd}
\\ & \,^2{\bf R}  \arrow{ru}\arrow[dashed]{rd}   & &  \arrow{rd} \bf C \\
 {\bf R}  \arrow{ru}  &&\arrow[dashed]{rd} && {\bf H} \arrow{ld}
 \\ &  {\bf C}   \arrow{lu}   & &   ^2{\bf H} \arrow{ld}  \\
 	&& \bf H \arrow{lu}  
\end{tikzcd} }
\end{table}

The Clifford clock contains in coded form exactly the same information about the matrix representation of $\G_{p,q}$ as does the classification table. Starting from the top $\bf R$ of each table, to get to any other vertex in the Classification Table, and the correponding ``Clifford Time", one can take steps to the right-down or to the left-down in the Classification Table, corresponding to advancing the same number of Clifford hours clockwise or counter-clockwise, respectively, on the clock. For the geometric algebra $\G_{p,q}$,  one proceeds from the top $\bf R$ to the right-down (clockwise) $0\le q \le 8$ steps (hours), followed by $0 \le p \le 8$ steps (hours) to the left-down (counterclockwish). The total number of steps taken (hours elapsed) determines $n=p+q$, the dimension of the matrix algebra over the base algebras (hours) on the Clifford Clock.

For example, we can use the Clifford Clock to construct an $8^{th}$ row of the Classifcation Table. Starting at Midnight, we count $8$ steps counter-clockwise, returning to $\bf R$. Next, again starting at Midnight, we take one step clock-wise, followed by $7$ steps counter-clockwise, landing at $\bH$. Next, two steps clock-wise, followed by $6$ steps counter-clockwise, again landing at $\bH$. We continue until the last entry is obtained by taking $8$ steps clockwise and again landing at $\bf R$. We record the $8^{th}$ row thus attained, together with the isomorphic geometric algebra $\G_{p,q}$ that it represents in Table 3. 
 \begin{table}
 	\centering
 	
 	\caption{Eighth row of geometric algebras $\G_{p,q}$.}
 	
 	\centering
 	\label{tab:table3}
 	\begin{tabular}{ccccccccc} \\
 		$\G_{8,0}$ & $\G_{7,1}$ &$ \G_{6,2}$  & $ \G_{5,3} $
 		& $ \G_{4,4} $   & $ \G_{3,5}$ & $ \G_{2,6}$ & $\G_{1,7}$ & $\G_{0,8}$ \\ 
 		$ M_{16}{\bf R} $  & $M_8\bH $& $M_8\bH$ & $ M_{16} {\bf R}$ & $M_{16} {\bf R}$ & $M_8\bH$ & $M_8\bH $& $M_{16}{\bf R} $  & $ M_{16}{\bf R} $
 	\end{tabular}  
 \end{table}

Our Axiom extols extending the real number system $\R $ to include new anti-commuting square roots of $\pm 1$. It is reasonable to ask what is the result of extending the complex numbers $\C $ to include new anti-commuting square roots of $\pm 1$?
 In this case, the complex geometric algebra
 \[ \G_{p,q}(\C):= \C(\be_1, \cdots , \be_p,\bbf_1, \ldots, \bbf_q  )= 
   \C(\be_1, \cdots , \be_p,i \bbf_1, \ldots, i \bbf_q  )= \G_{p+q}(\C ) ,   \] 
 so the study of {\it complex geometric algebras} is reduced to studying the
 structure of $\G_n(\C )$ for $n \ge 1$. Equivalently, we can say the entries in each row of the Classification Table 1 become algebraically isomorphic when considered over $\C $ instead of over $\R $. The building blocks in the complex case become $\G_0(\C )=\C $, and $\G_1(\C )=\, ^2\C$. The Classification Table for complex geometric algebras
 $\G_n(\C)$ is given in Table 4.
 
   The classification of complex geometric algebras in Table 4 has periodicity 2 rather than the periodicity 8 of the real geometric algebras in Table 1. The only new complex matrix entries in Table 4 are those entries containing $^2\C =M_2^D(\C )$. Some authors use the alternative equivalent notation $^2M_n(\C):=M_n(^2\C)$.

  \begin{table}
  	\centering
  	
  	\caption{Classification of Complex Geometric Algebras.}
  	
  	\centering
  	\label{tab:table4}
  	\begin{tabular}{cccccccc} \\
  		$\G_{0}(\C)$ &$\G_{1}(\C)$  &$\G_{2}(\C)$   &$\G_{3}(\C)$ 
  		&$\G_{4}(\C)$    &$\G_{5}(\C)$  & $\G_{6}(\C)$  &$\G_{7}(\C)$    \\ 
  		$\C $  & $^2\C $& $M_2(\C)$ & $ M_2(^2\C ) $ & $M_4 (\C )$ & $M_4(^2\C)$ & $M_8(\C ) $& $M_8(^2\C) $ 
  	\end{tabular}  
  \end{table} 

\section{The general linear group}

It is well-known that the general linear group $GL(N,\R)$ contains the {\it classical groups} as subgroups, such as the rotation group $SO(N,\R)$, the real orthogonal group $O(N,\R)$ of the real quadratic form on the Euclidean space $\R^N$, and the symplectic group $Sp(N,\R)$. In the listing of the real geometric algebras $\G_{p,q}$ given in Table 1, we find the matrix algebras $M_{N}(\R)$ for $N=2^k$, scattered throughout the table. In particular, down the middle of the table are the matrix
algebras $M_2(\R), M_4(\R), M_8(\R)$, isomorphic to
the geometric algebras $\G_{1,1}, \G_{2,2}, \G_{3,3}$, respectively, for which cases
$p=q=k$. More generally, we see that for each {\it even row} in Table 1, $p+q=2k$, and that no real matrix algebra appears in the {\it odd rows} of the table when $p+q=2k+1$.

Of course, each group $GL(N,\R )$ consists of the non singular matrices in
the corresponding matrix algebra $M_{2^k}(\R )$. Elie Cartan established a periodicity 8 relationship for (Clifford) geometric algebras in 1908, \cite[p.324]{LP97}. In 1959, R. Bott proved the periodicity 8 of homotopy groups of rotation groups. Along a different line of investigation into composition formulas for quadratic forms, A. Hurwitz and J. Radon established a periodicity 8 relationship in what is now known as the Hurwitz-Radon function.

Any positive integer can be expressed uniquely in the form $N=2^k k_0$, for non negative integers $k, k_0$ with $k_0$ odd.
The {\it Hurwitz-Radon function} $\rho(N):=\rho(2^k)$ is completely determined by the
diatic part $2^k$ of $N$, and the conditions
\beq \rho(N)=\pmatrix{2k+1 & {\rm if} & k=0 (mod \ 4) \cr 
	2k & {\rm if} & k = 1,2 (mod \ 4) \cr
	2k+2 & {\rm if}  & k=3 (mod \ 4)} \ \ \iff \ \ 
\rho(2^4 2^k)=\rho(2^k)+8. \label{HRnumbers} \eeq
The function $\rho(N)$ has the peculiar property that $\rho(N)=N$ only for $N=1,2,4,8$, which    crucial significance is discussed in the next section.  

In \cite{linkai2013}, Lin Kai-Liang considers the existence of anti-commuting matrices with squares $\pm 1$ in $GL(N, \bF )$, $O(N,\bF)$ and $Sp(N,\bF)$ where $\bF$ is any field of characteristic $\ne 2$, which characterize the Hurwitz-Radon function. He defines two
functions $R(N)$ and $S(N)$, closely related to the Hurwitz-Radon function $\rho(N)$ given in (\ref{HRnumbers}),
 \beq R(N):= \rho(N)-1, \ \ {\rm and} \ \ S(N):=\pmatrix{2k+1, & k=0 (mod\ 4) \cr
 	         2k, & k=1 (mod\ 4) \cr
 	         2k-1, & k=2,3(mod\ 4)}. \label{RSdef} \eeq
In particular, for the field $\bF=\R $, his Theorems 2.6 and 2.7 say there exists $r$ anti-commuting matrices in $GL(N,\R)$ with square $-1$, and $s$ anti-commuting matrices in $GL(N,\R)$ with square $+1$, iff $r \le R(N)$, and $s \le S(N)$.\footnote{Kai-Liang also gives a slightly more general result, his Theorem 2.5, which is closely related to our Theorem 2.} We now explain how Lin Kai-Liang's functions $R(N)$ and $S(N)$ are completely determined in Table \ref{tab:table5}, together with the Cartan periodicity 8 satisfied by the geometric algebras (\ref{cartan-period8}), and the periodicity 8 of $\rho(N)$, $R(N)$ and $S(N)$. We note also that $S(N)$ satisfies the useful
recurrence relation $S(2N)=R(N)+2$. 

\begin{table}
	\centering
	\caption{Geometric algebras isomorphic to $M_{2^k}(\R)$.}
	\label{tab:table5}	
	\[ \begin{tabular*}{\columnwidth}[t]{@{}c@{}c@{}c@{}c@{}c@{}c@{}c@{}c@{}c@{}c@{}c@{}c@{}@{}c@{}c@{}c@{}}
	$2k \ \  M_{2^k}(\R )  \quad \quad$  &  & &  &&$\G_{k+1,k-1} \quad $& & $\G_{k,k} \ \  $  & & $\G_{k-3,k+3}$ &\ \ $(+,-)$  & \ \ $S(N)$ \  &  &\ $R(N)$  & \\
	$ 0 \quad  M_{2^0}(\R )  \quad \quad$    && & &  &   & & $\G_{0,0}$ & & & $\ \ (1,0)$ & \ \ $1$& & $ 0 $ &  \\
	2$\quad M_{2^1}(\R ) \quad \quad$  &&&&&$\G_{2,0}$ && $\G_{1,1}$ && 
	&\ \ $(2,1)$ & $\ \ 2 $ && $ 1 $ & \\
	4$\quad  M_{2^2}(\R ) \quad \quad$ &&& &&$\G_{3,1}$ && $\G_{2,2}$&&  &$\ \ (3,2)$  & \ \ $3$ && $3$   \\
	6$\quad  M_{2^3}(\R ) \quad \quad$ &   &&  &&$\G_{4,2}$&&$\G_{3,3}$ && $\G_{0,6}$ & $\ \ (4,6)$  & $ \ \ 5 $   && $7$ \\
	$\quad i^2 \quad \quad \  $ &  & & &  & $-$ &  & $+$ & & $-$ & & & & & \\
	\end{tabular*}  \] 
\end{table}
\begin{table}
	\centering
	
	\caption{Geometric algebras isomorphic to $M_{2^4}(\R)$.}
	
	\centering
	\label{tab:table6}
	\begin{tabular}{ccccccccc} \\
		& 	$\G_{8,0}$  & $ \G_{5,3} $
		& $ \G_{4,4} $   & $\G_{1,7}$ &  $\G_{0,8}$ & $(+,-)$ & $S(2^4) $& $R(2^4)$ \\ 
		$i^2$ & $+$ & $-$ & $+$ & $-$ & $+$ & $(8,8)$  & $9$ & $8$
	\end{tabular}  
\end{table} 

Each row in Table \ref{tab:table5}, with the exception of the last row, gives all of the universal geometric algebras in Table \ref{tab:table1} that are isomorphic to the real matrix algebra $M_{2^k}(\R)$ for $k\ge 0$. The last row records the signs of the square of the pseudoscalar element $i$, carried over from the last row of Table \ref{tab:table1}. 
Since the pseudoscalar element $i$ is a $p+q=2k$-vector for each of the geometric algebras $\G_{p,q}$ in Table \ref{tab:table5}, it follows that it will anti-commute with each vector $\bv \in \G_{p,q}$, and, in particular, the $p$ and $q$ standard basis vectors with square $\pm 1$. The column of Table \ref{tab:table5} labeled by $(+,-)$ simply records the maximum number of standard basis elements with square $\pm 1$, respectively, that occur in the geometric algebras $\G_{p,q}$ in a given row. For example, in the row when $2k=4$, we find $(3,2)$, and in the row $2k=6$, we find $(4,6)$.

In the rows in Table \ref{tab:table5} when $k=0$, $S(1)=1$ and $R(1)=0$, and when
$k=1$, $S(2^1)=2$ and $R(2^1)=1$, in agreement with the value of $(+,-)$ for the universal geometric algebras in those rows, respectively. However, the argument is slightly more complicated for the rows when $k=2$ and $k=3$. In the row when $k=2$, $(+,-)=(3,2)$, 
corresponding to three anti-commuting real matrices with square $+1$ and two anti-commuting square matrices with square $-1$. However, recalling the case (\ref{qcase}) of quaternions, we can multiply the two real matrices together to get a third matrix which anti-commutes with the first two, so we can up the $2$ in $(3,2)$ to get $R(2^2)=3$. In the row when $k=3$, the value $(+,-)=(4,6)$ for the universal geometric algebras in that row. We can up the value of $4$ in $(4,6)$ to $S(2^3)=5$ by noting that the product of four anti-commuting real matrices with square $+1$ has itself square $+1$, and is a real matrix which anti-commutes with each of its component matrices. We can also up the $6$ in $(4,6)$ to $R(2^3)=7$ by noting that the pseudoscalar element in $\G_{0,6}$ anti-commutes with all the vectors in $\G_{0,6}$.

Table \ref{tab:table6} is a continuation of Table \ref{tab:table5} for $k=4$ and verifies the correctness of the Hurwitz-Radon functions given in (\ref{HRnumbers}) and
(\ref{RSdef}) as applied to geometric algebras isomorphic to $M_{2^k}(\R)$.

\section{Why is Bott periodicity so special?}

Bott periodicity, which has deep roots in homology theory, began with the study 
of {\it composition formulas} $[r,s,n]$ of {\it quadratic forms}, 
\beq (x^2_1 + x^2_2 + \cdots + x^2_r  )  (y_1^2 + y^2_ 2 + \cdots + y_s^2 )
 = z_1^2 + z_2^2  + \cdots + z_n^2. \label{compformula} \eeq
The $x_i,y_j$ are indeterminants and each $z_k:=z_k(X,Y)$ is a {\it bilinear}
form in the column vectors
 \beq X=\pmatrix{x_1 & \ldots & x_r }^T \ \ {\rm and} \ \ Y =\pmatrix{y_1 & \ldots & y_s }^T \label{basicquad} \eeq
with the coefficients in a field $\bf F$ with characteristic $\ne  2$.
Writing $Z =(z_1 $ $\ \ldots \ z_n )^T$ as a column vector, allows us to express
(\ref{basicquad}) as the matrix equation 
\[ (X^TX)( Y^TY) = Z^T Z . \]

Since $Z=Z(X,Y)$ is linear in both $X$ and $Y$, $Z = AY$ where $A$ is an $n\times s$ matrix linear in $X$.  The previous equation can then be written
\[  (X^T X)(Y^T Y) = Y^T A^T A Y , \]
 where
\beq A^T A =X^T X I_s = (x^2_1 + x^2_2 + \cdots + x^2_r  )I_s \label{eqnAAT} \eeq
 for the $s\times s$ identity matrix $I_s$.
Since $A$ is linear in $X$,
\[ A=x_1 A_1 + \cdots +x_r A_r \quad {\rm and} \quad A^T=x_1 A^T_1 + \cdots +x_r A^T_r \]
 where
each $A_i$ is a constant $n\times s$ matrix over $\bf F$.

 Substituting in the expressions for
$A$ and $A^T$ on the left side of equation (\ref{eqnAAT}), leads to the {\it Hurwitz Matrix Equations} (HME)
\beq A^T_i A_i = I_s \quad {\rm and} \quad A_i^T A_j + A^T_j A_i = 0 \label{HME1}  \eeq
for $i \ne j$. 
The Hurwitz Matrix Equations have a solution if and only if there is a
composition formula (\ref{compformula}) over $\bf F$ of size $[r, s, n]$.
Hurwitz \cite{hurwitz1923} and Radon \cite{radon1922}, studied the composition formulas $[r,s,n]$ for $s=n$. In this case, by letting $B_i=A_1^{-1}A_i$, the HME equations simplify to
\beq B_i^2 = -I_n \quad {\rm and} \quad B_iB_j+B_j B_i =0 \label{HME2} \eeq
for $i\ne j$, $2\le i, j\le n $, and $B_1=I_n$. 

The {\it Hurwitz-Radon Theorem}, named after work by Hurwitz first published in 1923, and work by Radon in 1922, states that there exists a composition formula of size $[r,n,n]$ if and only if $r\le \rho(n)$. Whereas Hurwitz and Radon considered composition formulas only over the fields $\R$ and $\C$, their theorem is also valid over any finite  
field $\bf F$ with characteristic $\ne 2$. The fact that $\rho(N)=N$ only for $N=1,2,4,8$, is intimately connected to fact that the only fields (and more generally division rings) are the real numbers $\R$, the complex numbers $\C $, the quaternions $\H $ and the octonians $\bf O$ which have the respective real dimensions $1, 2, 4,8$, respectively. A beautiful series of 3 Lectures by Daniel Shapiro \cite{dship1999}, gives a concise history of these results, together with an accounting of the current state of affairs. A much more detailed accounting can be found in his book \cite{dshipbook}, of which an ebook version is available on his homepage.   


Not surprisingly the the HME equations (\ref{HME2}), and even more general anti-commutative relationships between matrices \cite{linkai2013}, and their close relationship to the {\it composition of quadratic forms}, lead mathematicians to consider these problems in the context of Clifford algebras. Indeed, the general result that any finite dimensional real division algebra (not assuming any associativity or commutativity) has dimension 1,2,4, or 8 was established independently by Kervaire (1958), and Milnor and Bott (1958) using sophisticated parts of algebraic
topology that later became part of topological K-theory. For a discussion and references to this difficult and very technical problem, see \cite[pp. 186-87]{port1995}.
I have reached the point where I can only describe some of the other deep consequences of the theory that we have only touched upon here, and which represent the strenuous efforts of dozens of research mathematicians over a time span of more than 100 years. Perhaps the ``hairy ball problem" is the most famous -- the fact that it is impossible to have a continuous tangent, non-vanishing, vector field at every point on the $2$-sphere
$S^2$. The general statement of this theorem, first proved by Frank Adams in 1962, states that the maximal number of vector fields on the sphere $S^{N-1}$, linearly independent at each point is $\rho (N)-1$, where $\rho(N)$ is the Hurwitz-Radon function \cite{adams1961}.\footnote{In a personal communication, D. Shapiro pointed out that Adams' work can be viewed as an extension
	or refinement of the 1, 2, 4, 8 Theorem.}   

The Hurwitz-Radon function is also closely connected to the so called Hopf fibrations, which also only exist in dimensions $1,2,4,8$. A history of this problem, and its deep connections to homology theory, Morse theory, and K-theory is given in \cite{o-v2016}. A treatment of the Hopf fibration on the 2-sphere $S^2$, carried out in the geometric algebras $\G_3$ and $\G_4$, is given in  \cite{sob2015}.

\section*{Acknowledgements} I am grateful to D. Shapiro for invaluable remarks about the earlier version of this work, leading to many improvements.
 I thank Universidad de Las Americas-Puebla for many years of support.


\begin{thebibliography}{}
\bibitem{adams1961} J.F. Adams, {\em Vector Fields on Spheres}, Annals of Mathematics, Vol. 75, No. 3, May 1962.          	
\bibitem{b/w1935} R. Brauer and H. Weyl, {\em Spinors in n Dimensions}, 
American Journal of Mathematics, Vol. 57, No.2, 1935, pp 425-449.       
         	
\bibitem{budtraut} P. Budinich, A. Trautman, {\em The Spinorial Chessboard}, Springer (1988).   

\bibitem{cartan1966} E. Cartan, {\em The Theory of Spinors}, Dover, 1966.       	

\bibitem{wkc} W.K. Clifford, 1845-1879, Wikipedia,
https://en.wikipedia.org/wiki
/William\_Kingdon\_Clifford 

\bibitem{TD67} T. Dantzig, {\em NUMBER: The Language of Science},
Fourth Edition, Free Press, 1967.
\bibitem{dss1992} R. Delanghe, F. Sommen, V. Soucek, {\em Clifford Algebra and Spinor-Valued Functions: A Function Theory for the Dirac Operator}, Kluwer 1992.
\bibitem{grass1862} H. Grassmann, {\em Extension Theory}, A co-publication of the AMS and the London Mathematical Society.
http://www.ams.org/bookstore-getitem/item=HMATH-19
\bibitem{wrh} W.R. Hamilton, 1805-1865, Wikipedia.
https://en.wikipedia.org/wiki
/William\_Rowan\_Hamilton 

\bibitem{HD02} T. F. Havel, J.L. Doran, {\em Geometric Algebra in Quantum Information Processing},
Contemporary Mathematics, ISBN-10: 0-8218-2140-7, Vol. 305, 2002.
\bibitem{sources1967} G. Holton, General Editor, {\em Sources of Quantum Mechanics}, Dover 1967. 
\bibitem{hurwitz1923} A. Hurwitz, {\em \"Uber die Komposition der quadratischen Formen}, Math. Ann. 88, 1–25. Reprinted in Math. Werke, Bd. 2, Birkhäuser, 1962–3, 641–666.

\bibitem{LP97} P. Lounesto, 
\newblock {\em Clifford Algebras and Spinors, 2nd Edition}.
\newblock {Cambridge University Press}, Cambridge, 2001.
\bibitem{linkai2013} Lin Kai-Liang, {\em Hurwitz-Radon's symplectic analogy and Hua's cyclic recurrence relation}, Electronic Journal o Linear Algebra, ELA Vol. 26, pp. 858-872, 2013. 

\bibitem{o-v2016} V. Ovsienko and S. Tabachnikov, {\em Hopf Fibration and Hurwitz-Radon Numbers}, The Mathematical Intelligencer, No. 34, 4, pp. 11-18, 2016.
\bibitem{paulimat} Pauli Matrices, Wikipedia.
https://en.wikipedia.org/wiki/Pauli\_matrices

\bibitem{port1995} I. Porteous, {\em Clifford Algebras and the Classical Groups}, Cambridge studies in advanced mathematics, Cambridge University Press 1995.
\bibitem{radon1922} J. Radon, {\em Lineare Scharen orthogonaler Matrizen}, Abh. Math. Sem. Univ. Hamburg 1, 1–14. Reprinted in Collected Works, vol. 1, Birkhäuser 1987, 307–320.
\bibitem{dship1999} D. Shapiro, {\it Products of Sums of Squares}, A Series of
3 Lectures for a mini-course at the Universidad de Talca, December 1999. The Lectures
can be found on his homepage
 https://people.math.osu.edu/shapiro.6/
 \bibitem{dshipbook} D. Shapiro, {\it Compositions of Quadratic Forms}, W. de Gruyter Verlag, 2000. 
\bibitem{S0} G. Sobczyk, {\it The missing spectral basis in algebra and number theory}, The American
Mathematical Monthly 108 April 2001, pp. 336-346.

\bibitem{S1} G. Sobczyk, Hyperbolic Number Plane, {\it The College Mathematics
Journal}, Vol. 26, No. 4, pp.268-280, September 1995.

\bibitem{SNF} G. Sobczyk, {\em New Foundations in Mathematics: The Geometric Concept of Number},
\newblock Birkh\"auser, New York 2013.
	
\bibitem{S08} G. Sobczyk, Geometric Matrix Algebra, {\em Linear Algebra and its Applications}, 429 (2008) 1163-1173.

\bibitem{sob2015} G. Sobczyk, {\em Geometric Spinors, Relativity and the Hopf
	Fibration},  

\end{thebibliography}
\end{document}